\documentclass[11pt,reqno]{amsart}

\usepackage{amssymb,amsmath}

\newcommand{\be}{\begin{equation}}
\newcommand{\ee}{\end{equation}}
\newcommand{\bea}{\begin{eqnarray}}
\newcommand{\eea}{\end{eqnarray}}
\newcommand{\card}{\#}
\renewcommand{\P}{{\bf P}}
\newcommand{\es}{\emptyset}
\newcommand{\half}{\frac{1}{2}}
\newcommand{\eps}{\epsilon}

\newcommand{\R}{{\mathbb R}}

\newcommand{\E}{{\bf E\,}}

\def\Var{{\bf Var}}
\def\Cov{{\bf Cov}}
\def\Dk{{\mathcal D}}

\newtheorem{theorem}{Theorem}[section]

\def\sgn{{\rm sgn}}
\def\la{\langle}
\def\ra{\rangle}
\def\wx{\la w,X \ra}
\def\wxx{\la w,x \ra}
\def\one{{\mathbf 1}}
\def\t{t}
\theoremstyle{definition}

\theoremstyle{remark}

\numberwithin{equation}{section}

\begin{document}

\thispagestyle{empty}

\title[Stability of Weighted Majority]
{{Noise Stability of Weighted Majority}}

\author{Yuval Peres}\thanks{Research
partially supported by NSF grants \#DMS-0104073 and \#DMS-0244479.}

\address{Yuval Peres,
Department of Statistics, University of California, Berkeley. \newline
 {\tt peres@stat.berkeley.edu}}


\begin{abstract} Benjamini, Kalai and Schramm (2001)
showed that weighted majority functions of $n$ independent unbiased
bits are uniformly stable under noise: 
when each bit is flipped with probability $\epsilon$,
the probability $p_\epsilon$ that the weighted majority 
changes is at most
$C\eps^{1/4}$. They asked what is the best possible
exponent that could replace $1/4$.
We prove that the answer
is $1/2$. The upper bound obtained for $p_\epsilon$
is within a factor of $\sqrt{\pi/2}+o(1)$ from the known lower bound
when $\eps \to 0$ and $n\eps \to \infty$.
\end{abstract}

\maketitle

\section{{\bf Introduction}} \label{sec:intro}
In their study of noise sensitivity and stability of Boolean
 functions, Benjamini, Kalai and Schramm \cite{BKS}
showed that weighted majority functions of $n$ independent unbiased
${\pm 1}$-valued variables are uniformly stable under noise:  \newline
when each variable is flipped with probability $\epsilon$,
the weighted majority changes with probability at most
$C\eps^{1/4}$. They asked what is the best possible
exponent that could replace $1/4$.
In this note we prove
that the answer
is $1/2$.
Denote $\sgn(u)=u/|u|$ for $u \ne 0$ and $\sgn(0)=0$,
and let $N_\eps: \R^n \to \R^n$ be the noise operator that flips 
each variable in its input independently with probability $\eps$.
Formally, given a random vector $X=(X_1,\ldots,X_n)$,
the random vector $N_\eps(X)$
is defined as $(\sigma_1 X_1,\ldots,\sigma_n X_n)$
where the i.i.d.\ random variables $\sigma_i$ are independent of
$X$ and take the values $1,-1$ with probabilities $1-\eps,\, \eps$
respectively. 
\begin{theorem} \label{main}
Let $X=(X_1,\ldots,X_n)$ be a random vector uniformly distributed over
$\{-1,1\}^n$.
Given nonzero weights $w_1,\ldots, w_n \in \R$ 
and a threshold $\t \in \R$, consider the weighted majority function
$f: \R^n \to \{-1,0,1\}$ defined by
\be \label{wmaj}
f(x)=\sgn\Bigl(\sum_{i=1}^n w_i x_i-\t \Bigr)
\ee
Then for $\eps \le 1/2$,    
\be \label{half}
p_\eps(n,w,\t)=\P\Bigl(f(X) \ne f(N_\eps (X))\Bigr) \le 2\eps^{1/2}.
\ee
Moreover, $p_\eps^*=\limsup_{n \to \infty} \sup_{w,\t} p_\eps(n,w,\t)$
satisfies 
\be \label{asymp}
\limsup_{\eps \to 0} \frac{p_\eps^*}{\sqrt\eps} \le \sqrt{2/\pi} \,.
\ee
\end{theorem}
In the statement of the theorem we opted for a simple formulation:
Our proof yields the following sharper, but more involved estimate:
\be \label{basict}
p_\eps(n,w,t) \le \frac{2}{m} \E|B_m-\frac{m}{2}|
+[1-(1-\eps)^n] \binom{n}{\lfloor n/2\rfloor} 2^{-n} \,,
\ee
where $m=\lfloor \eps^{-1} \rfloor$ and $B_m$ is a Binomial$(m,1/2)$ variable.

\smallskip

It easy to see, and classical~\cite{Sh,Gu}, that for simple majority
(when all weights are equal) we have
\be \label{maj}
\lim_{n \to \infty} 
\P(\sgn \sum_{i=1}^n X_i \ne  \sgn \sum_{i=1}^n (N_\eps X)_i)
=\frac{1}{\pi}\arccos(1-2\eps)=\frac{2}{\pi}\sqrt{\eps}+O(\eps^{3/2})\,.
\ee
For the reader's convenience we include a brief argument: 

Since 
$\Cov\Bigl(\sum_{i=1}^n X_i,\, \sum_{i=1}^n( N_\eps X)_i\Bigr) =n(1-2\eps)$,
the central limit theorem implies that as $n \to \infty$,
$$
\frac{1}{\sqrt n} \Bigl(\sum_{i=1}^n X_i,\sum_{i=1}^n( N_\eps X)_i\Bigr)
\Rightarrow (Z_1, Z_1^*) \mbox{ \rm \; in law, }
$$
where $Z_1,Z_1^*$ are standard normals with covariance $1-2\eps$.
We can write $Z_1^*=Z_1\cos \alpha -Z_2\sin\alpha$
where $Z_1,Z_2$ are i.i.d. standard normals
and  $\alpha \in(0,\pi)$ satisfies $\cos \alpha=1-2\eps$.
Rotating the random vector $(Z_1,Z_2)$ 
by the angle $\alpha$ 
yields a vector with first coordinate $Z_1^*$.
Since $(Z_1,Z_2)$ has a rotationally-symmetric law,
the rotation changes the sign of the first coordinate with probability
$\alpha/\pi$. This verifies the left-hand side of (\ref{maj});
the right-hand side follows from Taylor expansion of cosine.


Thus the estimate (\ref{half}) is sharp
(up to the value of the constant).
Moreover, the ratio between the upper bound in 
(\ref{asymp}) and the value for simple majority in (\ref{maj})
tends to $\sqrt{\pi/2}<1.26$ as $\eps \to 0$.
We remark that the stability result in theorem \ref{main}
is stronger than an assertion about stability of half-spaces, 
$\{x \,:\, \sum_i w_ix_i>\theta \}$,
because we consider the weighted majority
as taking three values, rather than two.

\section{Proof of Theorem \ref{main}}
Using symmetry of $X_i$, we may assume that $w_i>0$
for $i=1,\ldots,n$.
Let $\wx=\sum_{i=1}^n w_i X_i$. 
We first consider the threshold $\t=0$. 
Later, we will extend the argument to thresholds $ t \ne 0$.

We will need the following well-known fact from~\cite{Er}:
\be \label{spern}
\P\Bigl(\wx=0\Bigr)\le \binom{n}{\lfloor n/2\rfloor} 2^{-n}.
\ee
Indeed, the collection $\Dk(w)$ of sets $D \subset \{1,\ldots n\}$
such that $\sum_{i \in D} w_i=\sum_{k \notin D} w_k$ forms an anti-chain
with respect to inclusion, so Sperner's theorem
(see \cite{AS}, Ch.~11) implies that the cardinality of $\Dk(w)$
is at most  $\binom{n}{\lfloor n/2\rfloor}$.
Finally, observe that a vector $x \in \{-1,1\}^n$ satisfies
$\wxx=0$ iff $\{i: x_i=1\}$ is in  $\Dk(w)$.

\medskip

Let $m=\lfloor \eps^{-1} \rfloor$ and let $\tau$ be a random
variable taking the values $0,1,\ldots, m$, with
$ \P(\tau=j) =\eps$ for $j=1,\ldots,m$ and  $\P(\tau=0) =1-m\eps$.
We use a sequence $\tau_1,\tau_2,\ldots,\tau_n$ of i.i.d.\ random variables
with the same law as $\tau$, to partition $[n]=\{1,\ldots, n\}$
into $m+1$ random sets
\be \label{partit}
A_j=\Bigl\{i \in [n] \, : \tau_i=j \Bigr\} \mbox{ \rm \; for } 0 \le j \le m.
\ee
Denote $S_j=\sum_{i \in A_j} w_iX_i$ and
let $Y_1=\sum_{i \notin A_j} w_iX_i =\wx-S_1$.
Observe that
$Y_1-S_1$ has the same law, given $X$, as $\la w, N_\eps(X) \ra$.
 Therefore, 
\bea \label{setup} 
p_\eps(n,w,0) &=& \P\Bigl( \sgn \wx \ne \sgn \la w, N_\eps(X) \ra \Bigr)
\\[1ex] \nonumber
&=& 
\P\Bigl( \sgn (Y_1+S_1) \ne \sgn (Y_1-S_1) \Bigr) \,.
\eea
Denote $\xi_j=\sgn(S_j)$. A key step in the proof
is the pointwise identity
\bea \label{keypoint}
 & &\one_{\displaystyle \{\sgn (Y_1+S_1) \ne  \sgn (Y_1-S_1) \}} 
\\[1ex]\nonumber
& & \; \; \; \; \; \; =2\cdot\one_{\displaystyle \{S_1 \ne 0 \}} \, \E
\Bigl(\half -\one_{ \displaystyle \{\sgn(S_1+Y_1) = -\xi_1 \}} 
\Big| \, Y_1, |S_1|\Bigr)\,.
\eea
To verify this, we consider three cases: \newline
\noindent{\bf (i)} Clearly both sides vanish if $S_1=0$. \newline
\noindent{\bf (ii)} Suppose that $0<|S_1| < |Y_1|$ and therefore 
$\sgn(Y_1+S_1)=\sgn(Y_1)$.
 The conditional distribution of $S_1$ given $Y_1$ and $|S_1|$
is uniform over $\{-|S_1|, |S_1|\}$, 
whence the 
conditional probability that $\sgn(S_1+Y_1) = -\xi_1$
is $1/2$. Thus both sides of (\ref{keypoint}) 
also vanish in this case.\newline
\noindent{\bf (iii)} Finally, suppose that $S_1 \ne 0$ and $|S_1| \ge |Y_1|$.
In this case $\sgn (S_1+Y_1) \ne -\xi_1$, so both sides of
(\ref{keypoint}) equal 1.

\medskip

Taking expectations in (\ref{keypoint}) 
and using (\ref{setup}), we deduce that
\bea \nonumber
p_\eps(n,w,0)&=& 2\,\E \Bigl[\one_{\displaystyle \{S_1 \ne 0 \}} 
\Bigl(\half -\one_{ \displaystyle \{\sgn \wx = -\xi_1 \} }\Bigr)\Bigr] 
  \\[1ex]
 &=& \frac{2}{m}\E \sum_{j \in \Lambda} 
\Bigl(\half -\one_{\displaystyle \{\sgn \wx=-\xi_j\}}\Bigr) \,,  \label{tight}
\eea
where $\Lambda=\{j \in [1,m]  : \, S_j \ne 0\}$.

The random variable $B_\Lambda=\card\{j \in \Lambda : \xi_j=1\}$
has a Binomial$(\card\Lambda, \half)$ distribution given $\Lambda$,
and satisfies the pointwise inequality
$$ 
\sum_{j \in \Lambda} 
\Bigl(\half -\one_{\displaystyle \{\sgn \wx=-\xi_j\}}\Bigr) \le 
\Bigl|B_\Lambda-\frac{\card \Lambda}{2}\Bigr|+  
\half \one_{\displaystyle \{\wx =0\}}
\sum_{j=1}^m \one_{\displaystyle \{A_j \ne \es\}}\,.
$$
To see this, consider the three possibilities for $\sgn \wx$.
Taking expectations and using (\ref{tight}), we get
\be \label{prep} 
p_\eps(n,w,0) \le  \frac{2}{m}
\E\Bigl|B_\Lambda-\frac{\card \Lambda}{2}\Bigr|+\P(A_1 \ne \es)\P(\wx=0).
\ee
Let $B_\ell$ denote a Binomial$(\ell,\half)$ random variable.
Since for any martingale $\{M_\ell\}_{\ell \ge 1}$
the absolute values $|M_\ell|$  
form a submartingale,
the expression $\E|B_\ell-\frac{\ell}{2}|$ is increasing in
$\ell $.
By averaging over $\Lambda$, we see that 
$\E|B_\Lambda-\frac{\card \Lambda}{2}| \le  \E|B_m-\frac{m}{2}|$.
In conjunction with (\ref{prep}) and (\ref{spern}), this implies
\be \label{basic}
p_\eps(n,w,0) \le \frac{2}{m} \E|B_m-\frac{m}{2}|
+[1-(1-\eps)^n] \binom{n}{\lfloor n/2\rfloor} 2^{-n}.
\ee

\medskip

Next, suppose that $f(x)=\sgn\Bigl(\sum_{i=1}^n w_i x_i-\t \Bigr)$,
where $t \ne 0$ is a given threshold.  
Let $X_{n+1}$ be a
$\, \pm 1 \, $ valued symmetric random variable, independent of 
$X=(X_1, \ldots, X_n)$, and define $w_{n+1}=t$.
Then 
\bea \label{halft}
p_\eps(n,w,\t) &=& \P\Bigl(f(X) \ne f(N_\eps (X))\Bigr) \\[1ex] \nonumber
&=&\P\Bigl(\sgn \sum_{i=1}^{n+1} w_iX_i \ne 
  \sgn  \bigl[\sum_{i=1}^{n}w_i(N_\eps X)_i+w_{n+1}X_{N+1}\bigr]\Bigr) \,,
\eea
and the argument used above to establish the bound 
(\ref{basic}) for $p_\eps(n,w,0)$, yields the same
bound for $p_\eps(n,w,\t)$. This proves (\ref{basict}).

To derive (\ref{half}), we may assume that $\eps \le 1/4$.
Use Cauchy-Schwarz to write
$ \E|B_m-\frac{m}{2}| \le \sqrt{\Var(B_m)}=\sqrt{m/4}$
and apply the elementary inequalities
$$\binom{n}{\lfloor n/2\rfloor} 2^{-n} \le \sqrt{3/4}\, n^{-1/2}, $$
(see, e.g., \cite{Pi}, Section 2.3) 
and $[1-(1-\eps)^n] \le \min\{n\eps,1\} \le \sqrt{n\eps}$,
to obtain
\be \label{now}
p_\eps(n,w,\t) \le m^{-1/2}+   \sqrt{n\eps} \cdot \sqrt{3/4}\, n^{-1/2}  \,.
\ee
Since $m=\lfloor \eps^{-1} \rfloor \ge 4/(5\eps)$ for $\eps \le 1/4$,
we conclude that
$$
p_\eps(n,w,\t) \le \Bigl(\sqrt{5/4}+ \sqrt{3/4}\Bigr) \, \eps^{1/2} < 
2\eps^{1/2} \,,
$$
and this proves (\ref{half}).

Finally, the central limit theorem implies that
 $$
\lim_{m \to \infty}
 \frac{\E|2B_m-m|}{\sqrt{m}}
 = \frac{1}{\sqrt{2\pi}}\int_{-\infty}^{\infty}
|u| e^{-u^2/2} \, du =  \sqrt{2/\pi} \,.
$$
This proves (\ref{asymp}).
\qed

\medskip

\noindent{\bf Remarks.}

\noindent{\bf 1.} Our proof  of Theorem \ref{main} was found in 1999, 
and was mentioned in \cite{BKS}. We present it here, with
more attention to the constants, in view of the recent interest
in related ``converse'' inequalities, see \cite{KKMO}.
The randomization idea which is crucial to the proof was inspired by
an argument of Matthews \cite{Mat} to bound cover times for Markov chains.
See also \cite{SW} for related random walk estimates.

\noindent{\bf 2.} After I described the proof  of Theorem \ref{main}
to R. O'Donnell, he found (jointly with A. Klivans and R. Servedio)
some extensions and applications of the argument 
to learning theory, see \cite{KOS} for this and many other results.

\noindent{\bf 3.}  The proof of Theorem \ref{main} extends verbatim to the case
 where $X_i$ are independent symmetric real-valued random variables
with $\P(X_i=0)=0$ for all $i$. However, this extension reduces to  
Theorem \ref{main} by conditioning on $|X_i|$. A more 
interesting extension would 
be to replace the symmetry assumption on $X_i$ by the assumption $\E X_i$=0.

\noindent{\bf 4.} Is simple majority the most noise sensitive of the weighted 
majority functions (asymptotically when $\eps \to 0$ and $n\eps \to \infty$) ?
\newline
In particular, is it possible to replace the right-hand side
of (\ref{asymp}) by $2/\pi$? 


\noindent{\bf Acknowledgement.} I am grateful to
I. Benjamini, G. Kalai and O. Schramm for suggesting the problem,
and to E. Mossel, R. Peled, O. Schramm,
R. Siegmund-Schultze and H. V. Weizs\"acker for useful discussions.

\bibliographystyle{amsplain}

\end{document}